\documentclass[10pt,leqno]{jaums}
\usepackage[latin1]{inputenc}
\usepackage{amsmath,latexsym,amssymb,graphicx,dsfont}
\usepackage{natbib}

\theoremstyle{thmit} 
\newtheorem{Th}{Theorem}[section]
\newtheorem{Lemma}[Th]{Lemma}
\newtheorem{Cor}[Th]{Corollary}
\newtheorem{Prop}[Th]{Proposition}

\theoremstyle{thmrm} 

\newtheorem*{oldproof}{Proof}
\renewenvironment{proof}[1][{}]{\begin{oldproof}[#1]}{\qed\end{oldproof}}

\title{Rate of Escape of Random Walks on Free Products}
\author{Lorenz A. Gilch}
\address{University of Technology Graz, Institut für Mathematik C, Steyrergasse
  30,\\ A-8010 Graz, Austria}
\email{gilch@TUGraz.at}
\url{http://www.math.tugraz.at/$\sim$gilch/}


\begin{document}

\maketitle

\begin{abstract}
Suppose we are given the free product $V$ of a finite family of finite or
  countable sets $(V_i)_{i\in\mathcal{I}}$
  and probability measures on each $V_i$, which govern random walks on it. We
  consider a transient random walk on the free product arising naturally from
  the random walks on the $V_i$. We prove the existence of the rate of
  escape with respect to the block length, that is, the speed, at which the
  random walk escapes to infinity, and furthermore we compute formulas for
  it. For this purpose, we present three different techniques providing three
  different, equivalent formulas. 
\end{abstract}

\section{Introduction}

Consider a transient Markov chain $(Z_n)_{n\in\mathbb{N}_0}$ on a state space $V$ and a
suitable length function $l$ on $V$ representing a `word length' with respect
to the
starting point of the Markov chain. We are interested in whether
the sequence of random variables $l(Z_n)/n$ converges almost surely to a
constant, and if so, to compute this constant. If the limit exists, it is called the
\textit{rate of escape}, or the \textit{drift with respect to $l$}. In this paper, we study this
question for random walks on general free products.
\par
To outline some background material, on the $d$-dimensional grid $\mathbb{Z}^d$, where \mbox{$d\geq 1$}, random
walks can be described by the sum of $n$
independent and identically distributed random variables, the increments of $n$ steps. By the weak law of
large numbers the limit $\lim_{n\to\infty} |Z_n|/n$, where $|\cdot|$ is the
distance on the grid to the starting point of the random walk, exists almost
surely. Furthermore, this limit is positive if the increments have non-zero
mean vector. 
\par
It is well-known that the rate of escape exists also for transitive random walks on finitely
generated groups, where the random walks arise from a
probability measure on the group elements. This follows from \textit{Kingman's subadditive
  ergodic  theorem}; see Kingman \cite{kingman}, Derriennic \cite{derrienic} and
Guivarc'h \cite{guivarch}. If $l$ is the metric of the Cayley graph, then the
limit $\lim_{n\to\infty} l(Z_n)/n$ exists almost surely and is positive. There
are many detailed results for random walks on groups and wreath products: Mairesse
\cite{mairesse-mantaray} computed a explicit formula in terms of the unique
solution of a system of polynomial equations for the rate of escape of
 random walks on the braid group. Lyons,
Pemantle and Peres \cite{lyons-pemantle-peres} gave a lower bound for the rate
of escape of inward-biased random walks on lamplighter groups. Dyubina
\cite{dyubina} proved that the drift on the wreath product $A\wr \mathbb{Z}/2$ is
zero, where $A$ is a finitely generated group. Erschler \cite{erschler2}
investigated asymptotics of the drift of symmetric random walks on finitely
generated groups. An important link
between drifts and harmonic analysis was obtained by Varopoulos
\cite{varopoulos}. He proved that for symmetric finite range random walks on groups the existence of non-trivial bounded
harmonic functions is equivalent to a non-zero rate of escape. This leads to a
link between the rate of escape and the entropy of random walks, compare
e.g. with Kaimanovich and Vershik \cite{kaimanovich-vershik} and Erschler
\cite{erschler2}. The rate of escape has also been studied on trees: Cartwright, Kaimanovich and Woess \cite{cartwright-kaimanovich-woess} investigated the
boundary of homogeneous trees and the drift on them. Nagnibeda and Woess \cite[Section 5]{woess2} proved that the
rate of escape of random walks on trees with finitely many cone types is non-zero and give a
formula for it.
\par
For a restricted class of free
products of finite groups, Mairesse \cite{mairesse} and Mairesse and Math\'eus
\cite{mairesse1} have developed a specific technique for computation of
the above limit with respect to the word length. These papers were the starting point for the present investigation of
arbitrary free products. We consider the free product of finitely many sets, on
which Markov chains are given, and construct in a natural way a random walk on the
free product. The techniques we use for rewriting probability generating functions in
terms of functions on the factors of the free product were introduced independently and
simultaneously by Cartwright and Soardi \cite{cartwright-soardi}, Woess \cite{woess3}, Voiculescu \cite{voiculescu} and McLaughlin \cite{mclaughlin}.
\par
Our aim is to show the existence of the above rate of escape
$\ell$ with respect to the word
length, and also to compute formulas for it.
For this purpose, we will present three different,
equivalent formulas for $\ell$  using three different techniques. In Section \ref{bl_fet} we
prove existence and a formula for $\ell$ by purely probabilistic reasoning. In Section \ref{bl_dgf} we compute the proposed limit using
double generating functions and applying a theorem of Sawyer and Steger \cite[Theorem 2.2]{sawyer}. The third
approach for the computation of $\ell$ in Section \ref{bl_lp} works only for
free products of finitely generated groups and is based on a technique
which was already used by Ledrappier \cite{ledrappier} and Furstenberg \cite{furstenberg}. Section \ref{samples} presents
sample computations and in Section \ref{remarks} we give
additional remarks about extensions of these techniques to further results.

\section{Free Products}
\subsection{Free Products and Random Walks}
Let $\mathcal{I}:=\{1,\dots,r\}$, $r\geq 2$.
Consider $r$ random walks with
transition matrices $P_i$
on pairwise disjoint finite
or countable state spaces $V_i$, where $i\in\mathcal{I}$. The corresponding single and $n$-step transition
probabilities are denoted by $p_i(x,y)$ and $p_i^{(n)}(x,y)$, where
$x,y\in V_i$. For every $i\in\mathcal{I}$ we select an element $o_i$ of $V_i$
as the `root'. To help visualize this, we
think of graphs $X_i$ with vertex sets $V_i$ and roots $o_i$ such that there is
an oriented edge $x\to y$ if and only if $p_i(x,y)>0$. Furthermore, we shall
assume that for every $i\in\mathcal{I}$ and every $x\in V_i$ there is an
$n\in\mathbb{N}$ such that $p^{(n)}_i(o_i,x)>0$. For the sake of simplicity we
assume $p_i(x,x)=0$ for every $i\in\mathcal{I}$ and $x\in V_i$. 
\par
Let $V_i^\times := V_i\setminus\{o_i\}$ for every $i\in\mathcal{I}$. The next
step is the construction of a new Markov chain on the \textit{free product} 
$V:=V_1\ast \dots \ast V_r$, the set of `words'
\begin{equation}\label{freeproduct}
x=x_1x_2\dots x_n
\end{equation}
with letters, also called \textit{blocks}, from the sets $V_i^\times$
such that no two successive letters come from the same
  $V_i$. The empty word $o$ describes the root of $V$. If $u=u_1\dots u_m\in V$ and $v=v_1\dots
v_n\in V$ with $u_m\in V_i$ and $v_1\notin V_i$, then $uv$ stands for their
concatenation as words. We also define $vo=ov=v$ for all $v\in V$. We regard
each $V_i$ as a subset of $V$, identifying each $o_i$ with $o$.
\par
We lift $P_i$ to a transition matrix $\bar P_i$ on $V$: if
$z\in V$ is $o$ or has last letter not in $V_i$, and if $v,w\in V_i$, then we
set $\bar
p_i(zv,zw):=p_i(v,w)$. Otherwise we set $\bar p_i(x,y):=0$. We choose
$0<\alpha_1,\dots,\alpha_r\in\mathbb{R}$ with $\sum_{i\in\mathcal{I}} \alpha_i =
1$. Then we obtain a new transition matrix on $V$ given by
$$
P=\sum_{i\in\mathcal{I}} \alpha_i \bar P_i.
$$
The random walk governed by $P$ is described by the sequence of random variables
$(Z_n)_{n\in\mathbb{N}_0}$. The associated single and $n$-step transition
probabilities are denoted by $p(x,y)$ and $p^{(n)}(x,y)$ for $x,y\in V$.
\par
Let $x=x_1\dots x_m\in V\setminus\{o\}$. The \textsf{type} $\tau(x)$ of $x$ is
defined to be $i$ if $x_m\in V_i$. The \textsf{block length} $\ell(x)$ of $x$
is defined to be $m$. We also set $\tau(o):=0$ and $\ell(o):=0$. We want to show existence of the $P$-almost sure limit
$\ell=\lim_{n\to\infty} \ell(Z_n)/n$, the \textit{rate of escape
  with respect to the block length},
and to present formulas for it.  Let $k\in\mathbb{N}$. Then $x^{(k)} := x_1\dots
  x_j$, where $j=\min\{k,m\}$, is the \textsf{truncation} at length $k$. We also
  write $x^\bot := x_1\dots x_{m-1}$ for the truncation at length $\ell(x)-1$, if
  $x\neq o$. Furthermore we denote by $\tilde x := x_m$ the \textsf{terminal
    block} of $x$, if $m>0$, and set $\tilde o:=o$. The \textsf{cone} rooted at $x$ is the set 
$$
C_x := \bigl\lbrace y\in V \ |\  y^{(m)}=x\bigr\rbrace \subseteq V.
$$
If $y\in V_i$, then the set of
  \textsf{successors} of $y$ is given by
$$
\mathcal{S}(y)  :=  \bigl\lbrace w\in V_i \ |\ p_i(y,w)>0\bigr\rbrace
$$
and the set of \textsf{predecessors} by
$$
\mathcal{P}(y) := \bigl\lbrace w\in V_i \ |\ p_i(w,y)>0\bigr\rbrace.
$$
\par
We now introduce some probability generating functions. For this purpose, let 
\mbox{$T_y:=\min \{k\geq 0\,|\,Z_k=y\}$}, resp. \mbox{$S_y:=\min \{k>
  0\,|\,Z_k=y\}$}, be the stopping
time of the first visit, resp. the first return to $y\in V$. Denote by $\mathbb{P}_x$ the probability
measure on $V^{\mathbb{N}_0}$ that governs the random walk starting at $x\in V$. For
$z\in\mathbb{C}$ and $x,y\in V$ let
\begin{eqnarray*}
G(x,y|z)& = &\sum_{n\geq 0} p^{(n)}(x,y)\,z^n, \quad \quad F(x,y|z)=\sum_{n\geq 0}
\mathbb{P}_x[T_y=n]\,z^n,\\
 U(x,y|z)& =&\sum_{n\geq 1} \mathbb{P}_x[S_y=n]\,z^n,\\
 L(x,y|z)& = & \sum_{n\geq 0} \mathbb{P}_x\bigl[\forall k\in \{1,\dots,n-1\}: Z_k\neq x,Z_n=y\bigr]\,z^n.
\end{eqnarray*}
The analogous functions for the random walks on the single factors $V_i$
are denoted by $G_i(u,v|z)$, $F_i(u,v|z)$ and $L_i(u,v|z)$, where $u,v\in V_i$.  We make the basic
assumption that the radius of convergence of $G(o,o|z)$ is greater than $1$, which
implies transience of our random walk on $V$. Thus, we may exclude the case
$r=2=|V_1|=|V_2|$. This convergence property is
fulfilled if each $p_i$ is reversible and due to non-amenability also for random walks on free products of finitely
generated groups,
where the $P_i$ depend only on a probability measure on the 
single groups. Note that for $|z|<1$, $\sum_{y\in V_i} G_i(x,y|z)=1/(1-z)$ for
every $i\in\mathcal{I}$ and all $x\in V_i$. This will be used several times in
the sequel.
\par
Recall the following equations:
\begin{Lemma}\label{lemma-props}
Let $x,y\in V$, $w\in V\setminus\{o\}$ such that $w^{(1)}\notin V_{\tau(x)}$, and $z\in\mathbb{C}$. Then:
\begin{eqnarray*}
(i) & G(x,x|z)  &=   \frac{1}{1-U(x,x|z)}\,, \\
(ii) & G(x,y|z)  & =   F(x,y|z) \cdot G(y,y|z) \,,\\
(iii) &
G(x,y|z)  &=  G(x,x|z) \cdot L(x,y|z)\,,\\
(iv)& F(o,xw|z)&= F(o,x|z)\cdot F(x,xw|z)\,,\\
(v)& L(o,xw|z) &= L(o,x|z)\cdot L(x,xw|z)\,,\\
(vi)& L(x,xw|z) &= L(o,w|z).
\end{eqnarray*}
Equations analogous to $(i)$, $(ii)$, $(iii)$ hold for the generating functions on the single factors
$V_i$ for every $i\in\mathcal{I}$.
\end{Lemma}
\begin{proof} 
For $(i)$ and $(ii)$ see Woess \cite[Lemma 1.13]{woess}. Equation
$(iii)$ is obtained by conditioning with respect to the last visit at $x$
before finally walking to $y$. $(iv)$ and $(v)$ are obtained by conditioning with respect
to the first resp. last visit at $x$, which must be visited before finally
walking to $xw$. $(vi)$ holds, as due to the tree-like structure of the free product the probability of walking from $x$ to $xw$ in $n$ steps without returning
to $x$ is the same as walking from $o$ to $w$ in $n$ steps without returning to $o$.
\end{proof} 

We now explain the correspondence between $F(x,y|z)$ and $F_i(x,y|z)$, resp. $L(x,y|z)$ and $L_i(x,y|z)$. Therefore define for
$i\in\mathcal{I}$ and $z\in\mathbb{C}$
$$
\bar H_i(z)  :=  \sum_{n=2}^\infty \mathbb{P}_o[S_o=n,Z_1\notin V_i]\,z^n \
\textrm{ and } \ 
\xi_i(z) :=  \frac{\alpha_i z}{1-\bar H_i(z)}\,;
$$
see Woess \cite[Proposition 9.18]{woess}.
Note that $\bar H_i(1)$ is the probability of starting at some $x\in V_i$
and returning to the same $x$ without having visited a neighbour
$y\in\mathcal{S}(x)$ of $x$. Similarly $\xi_i(1)$ is the probability of
starting at some $x\in V_i$ and at some time visiting a neighbour
$y\in\mathcal{S}(x)$. Observe that for positive $z$ the functions $\bar H_i(z)$ and $\xi_i(z)$ are strictly increasing inside
their radii of convergence, which are greater than $1$.
\begin{Lemma}\label{lemma-fi}
Let $i\in\mathcal{I}$, $x,y\in V_i$ and $z\in\mathbb{C}$. Then
$$
(i)\ F(x,y|z)=F_i\bigl(x,y|\xi_i(z)\bigr) \quad \textrm{ and } \quad (ii)\
L(x,y|z)=L_i\bigl(x,y|\xi_i(z)\bigr).
$$
\end{Lemma}
For the proof of $(i)$ see Woess \cite[Proposition 9.18~(c)]{woess}. Statement
$(ii)$ is proved analogously.

\begin{Lemma}
\label{lemma-prop}
$\xi_i:=\xi_i(1)<1$ for all $i\in\mathcal{I}$.
\end{Lemma}
\begin{proof}
Let $H_i(z):=U(o,o|z)-\bar H_i(z)$. By transience we have
$$
U(o,o|1)=\sum_{i\in\mathcal{I}} H_i(1) < 1.
$$
Furthermore 
$$
H_i(1)=\alpha_i \sum_{s\in \mathcal{S}(o_i)} p_i(o_i,s)
\underbrace{F(s,o|1)}_{\leq 1} \leq \alpha_i.
$$
Hence,
$$
\xi_i = \frac{\alpha_i}{1-\sum_{j\in\mathcal{I}\setminus\{i\}} H_j(1)} \leq
\frac{\alpha_i}{1-\sum_{j\in\mathcal{I}\setminus\{i\}} \alpha_j} =
\frac{\alpha_i}{1-(1-\alpha_i)} = 1.
$$
Observe that $\xi_j<1$ for all $j\in\mathcal{I}\setminus\{i\}$, if
$H_i(1)<\alpha_i$ for some $i\in\mathcal{I}$. Assume $H_i(1)=\alpha_i$ for some
$i\in\mathcal{I}$. Then $\bar H_i(1) = U(o,o|1)-H_i(1) < 1-\alpha_i$,
and thus there is $j\in\mathcal{I}\setminus\{i\}$ such that $H_j(1)<
\alpha_j$. Thus $\xi_i<1$. Since
$$
H_i(1)=\alpha_i \sum_{s\in \mathcal{S}(o_i)} p_i(o_i,s)
F(s,o|1) = \alpha_i
$$
we have $F(s,o|1)=1$ for all $s\in\mathcal{S}(o_i)$. But now we obtain the contradiction
$$
F(s,o|1)=F_i(s,o_i| \xi_i)<1,
$$
as $\xi_i<1$ and $F_i(s,o_i|x)$, $0\leq x\in\mathbb{R}$, is strictly increasing with \mbox{$F_i(s,o_i|1)\leq 1$}. 
\end{proof}

\subsection{Limit of the Random Walk}
As we have assumed transience for the random walk on $V,$ the random walk escapes
to infinity in the sense that almost surely every finite $A\subseteq V$ is visited only
finitely often.
We shall now investigate the route of the escape of the random walk on $V$, which provides the main tool for further computations.
Define for $x\in V$, $i\in\mathcal{I}$ with
$\tau(x)\neq i$ and $S\subseteq V_i$ the set $xS :=\{ xy \ |\  y\in S\}$. Then we obtain:
\begin{Lemma}
\label{lemma-copyend}
If $x\in V$ and $\tau(x)\neq i$ then
$$
\mathbb{P}_o\bigl[ Z_n\in xV_i \ \textrm{holds for infinitely many}\ n \bigr] =
0.
$$
\end{Lemma}
\begin{proof}
By Lemmas \ref{lemma-props} and \ref{lemma-fi},
\begin{eqnarray*}
\sum_{n\geq 0} \mathbb{P}_o[Z_n\in xV_i] & = & \sum_{y\in V_i} G(o,xy|1) \\
& = & \sum_{y\in V_i} G(o,o|1)\, L(o,x|1)\, L_i(o_i,y|\xi_i) \\
& = & G(o,o|1)\cdot L(o,x|1) \cdot \sum_{y\in V_i}
\frac{G_i(o_i,y|\xi_i)}{G_i(o_i,o_i|\xi_i)} \\
& = & \frac{ G(o,o|1)\ L(o,x|1)}{G_i(o_i,o_i|\xi_i)} \frac{1}{1-\xi_i} < \infty.
\end{eqnarray*}
The Borel-Cantelli lemma implies the proposed statement. 
\end{proof} 
Now we are able to specify how the random walk on $V$ escapes to
infinity. Let $V_\infty$ denote the set of infinite words $x_1x_2\dots$ in
which each of the letters $x_j$ belongs to $\bigcup_{i\in\mathcal{I}}
V_i^\times$, no consecutive letters come from the same $V_i^\times$, and
infinitely many letters come from each $V_i^\times$. Then we obtain:
\begin{Prop}
$\ell(Z_n)$ tends to infinity $\mathbb{P}_o$-a.s., as
$n\to\infty$. Furthermore, there exists a $V_\infty$-valued random variable $Z_\infty$, such that
$$
\lim_{n\to\infty} Z_n = Z_\infty \quad \mathbb{P}_o-\textrm{a.s.,}
$$
with convergence in the sense that the length of the common prefix of $Z_n$
and $Z_\infty$ tends to infinity.
\label{elltoinfty}
\end{Prop}
\begin{proof}
We prove by induction that for each $m\in\mathbb{N}$ there is almost surely
some $n_m\in\mathbb{N}$ with $\ell(Z_{n_m})=m$ and $\ell(Z_n)> m$ for all $n>
n_m$. By Lemma \ref{lemma-copyend}, the random walk visits the state set
$\bigcup_{i\in\mathcal{I}} V_i$ finitely often $\mathbb{P}_o$-a.s.. Therefore
there is almost surely some $n_1\in\mathbb{N}$ such that  $Z_{n_1}\in\bigcup_{i\in\mathcal{I}} V_i$
  $\mathbb{P}_o$-a.s. and $Z_{n}\notin\bigcup_{i\in\mathcal{I}} V_i$ for all
  $n>n_1$. Thus $\ell(Z_{n_1})=1$ and $\ell(Z_n)> 1$ for all $n>n_1$. \\
Assume
  now that $\ell(Z_{n_m})=m$ and $\ell(Z_n)> m$ for all
  $n> n_m$ and let \linebreak[4]\mbox{$\mathcal{I}':=\mathcal{I}\setminus \{\tau(Z_{n_m})\}$}.
Again by Lemma \ref{lemma-copyend} the random walk visits the
  state set $\bigcup_{i\in\mathcal{I}'} Z_{n_m}V_i$
  finitely often $\mathbb{P}_o$-a.s.. Then there is almost surely some
  $n_{m+1}\in\mathbb{N}$ such that
  $Z_{n_{m+1}}\in\bigcup_{i\in\mathcal{I}'} Z_{n_m}V_i$ and $Z_{n}\notin\bigcup_{i\in\mathcal{I}'} Z_{n_m}V_i$ for all $n>n_{m+1}$. Thus $\ell(Z_{n_{m+1}})=m+1$ and $\ell(Z_n)>
  m+1$ for all $n> n_{m+1}$. Thus $\ell(Z_n)$ tends to infinity, as $n\to\infty$.\\
Obviously the sequence $(Z_n)_{n\in\mathbb{N}_0}$ converges to an infinite
word in $V_\infty$ with $Z_\infty^{(m)}=Z_{n_m}$ for all
$m\in\mathbb{N}$. \\
\end{proof}

%
%

\section{Exit Time Technique}
\label{bl_fet}
In this section we investigate the random walk on $V$ in detail, prove the
existence of $\ell$ and derive a formula for it. The following technique was motivated by Nagnibeda and Woess \cite[Section 5]{woess2}. Let $k\in\mathbb{N}$. The \textsf{exit time with respect to the block length $k$} is
\begin{displaymath}
\mathbf{e}_k := \min \Bigl\lbrace m\in\mathbb{N}\ \Bigl|\  \forall n\geq m: Z_n^{(k)} \textrm{
  constant }\Bigr\rbrace.
\end{displaymath}
In particular $\mathbf{e}_0=0$. The \textsf{exit point with respect to the block length $k$}
is \mbox{$W_k := Z_{\mathbf{e}_k}$}.  Thus, $\mathbf{e}_k$ is the first instant from
which point the first $k$ blocks remain constant, and $W_k=x$ if and only if
at time  $\mathbf{e}_k-1$ the random walk is at state $x_1\dots x_{k-1}s$ with
some $s\in \mathcal{P}(x_k)$, at time $\mathbf{e}_k$ at state $x$, and
thereafter remains in the \mbox{cone $C_x$}.
\par
As $Z_n$ converges almost surely to a random variable $Z_\infty$ with values in
$V_\infty$, we have $\mathbf{e}_k\to\infty$ as $k\to\infty$, almost surely. The
\textsf{$k$-th increment} is $\mathbf{i}_k := \mathbf{e}_k-
\mathbf{e}_{k-1}$ and the \textsf{maximal temporary
  exit time at time $n$} is defined as
\begin{displaymath} 
\mathbf{k}(n):= \max \bigl\lbrace k\in\mathbb{N}_0\ |\ \mathbf{e}_k\leq n\bigr\rbrace.
\end{displaymath}
Define now for $i\in\mathcal{I}$, $y\in V$ and $n\in\mathbb{N}_0$
$$
k_i^{(n)}(o,y) := \mathbb{P}_o\bigl[\forall j\in\{0,\dots,n\}: Z_j\notin V_i^\times,Z_n=y\bigr]
$$
and the corresponding generating function
\begin{equation}\label{ki}
K_i(o,y|z)  :=  \sum_{n=0}^\infty k_i^{(n)}(o,y)\, z^n =
\sum_{n\geq 0} \bar H_i(z)^n \cdot L(o,y|z) = \frac{L(o,y|z)}{1-\bar H_i(z)}\,.
\end{equation}
Now we can prove:
\begin{Prop}
$(W_k,\mathbf{i}_k)_{k\in\mathbb{N}}$ is a Markov chain with transition probabilities
\begin{eqnarray*}
&& \mathbb{P}_o\bigl[W_{k+1}=w_{k+1},\mathbf{i}_{k+1}=n_{k+1}\, \bigl|\,
W_k=w_k,\mathbf{i}_k=n_k\bigr]  \\[2ex]
& = & \frac{1-\xi_{\tau(w_{k+1})}}{1-\xi_{\tau(w_k)}} \cdot
\sum_{s\in \mathcal{P}(y)} 
\Bigl[ k_{\tau(w_k)}^{(n_{k+1}-1)}(o,s)\cdot p(s,y) \Bigr]
\end{eqnarray*}
for $n_k,n_{k+1}\in\mathbb{N}$, $w_k=x_1\dots x_k\in V$, $w_{k+1}=w_ky$, where $y\in \bigcup_{i\in\mathcal{I}\setminus\{\tau(w_k)\}}V_i^\times$.  
\end{Prop}

\begin{proof}
Define $\bar V_i=\bigcup_{j\in\mathcal{I}\setminus\{i\}} V_j^\times$. Let
  $w_0=o$, $w_1= g_1 \in \bigcup_{i\in\mathcal{I}}V_i^\times$ and \mbox{$w_i=w_{i-1}g_i$} with
  $g_i\in \bar V_{\tau(w_{i-1})}$ for $2\leq i \leq k$.

For $i\in\{1,\dots,k\}$ the inclusion $[W_{i+1}=w_{i+1}]\subseteq [W_i=w_i]$
holds, as $w_{i+1}$ determines the element $w_i$ uniquely. 
Let
$n_1,\dots,n_{k+1}\in\mathbb{N}$ and write for \mbox{$m\in\{k,k+1\}$}
\begin{displaymath}
\bigl[W_1^{m}=w_1^m,\mathbf{i}_1^{m}=n_1^m\bigr] := 
\bigl[\forall j\in\{1,\dots,m\}: W_j=w_j,\mathbf{i}_j=n_j\bigr]. 
\end{displaymath}
This event can be described as follows: start at $o$, walk in $n_1-1$ steps to
a predecessor of $w_1$ inside $V_{\tau(w_1)}$, then walk to $w_1$; then stay
inside $C_{w_1}$ and walk in $n_2-1$ steps to a vertex in
$w_1\mathcal{P}(g_2)$, from there to $w_2$, and so on. More formally we obtain,
writing $n_1^s=\sum_{t=1}^s n_t$:
\begin{eqnarray*}
& & \mathbb{P}_o[W_1^k=w_1^k,\mathbf{i}_1^k=n_1^k] \\
&=&
\mathbb{P}_o \left[
\begin{array}{c}
\forall \lambda\in\{0,\dots,k-1\}\forall j\in\{1,\dots,n_\lambda
-1\}:\\
Z_{n_1^\lambda+j}\in C_{w_\lambda},Z_{n_1^{\lambda+1}-1}\in w_{\lambda}\mathcal{P}(g_{\lambda+1}),Z_{n_1^{\lambda+1}}=w_{\lambda+1}
\end{array} \right] \cdot \bigl(1-\xi_{\tau(w_k)}\bigr).
\end{eqnarray*}
Analogously,
\begin{eqnarray*}
&& \mathbb{P}_o[W_1^{k+1}=w_1^{k+1},\mathbf{i}_1^{k+1}=n_1^{k+1}]\\
&=&\mathbb{P}_o \left[
\begin{array}{c}
\forall \lambda\in\{0,\dots,k-1\}\forall j\in\{1,\dots,n_\lambda
-1\}:\\
Z_{n_1^\lambda+j}\in C_{w_\lambda},Z_{n_1^{\lambda+1}-1}\in w_{\lambda}\mathcal{P}(g_{\lambda+1}),Z_{n_1^{\lambda+1}}=w_{\lambda+1}
\end{array} \right] \cdot\\
&& \quad
 \mathbb{P}_{w_k}\left[
\begin{array}{c}
\forall j\in \{1,\dots,n_{k+1}-2\}:Z_j\in C_{w_k},\\
Z_{n_{k+1}-1}\in
w_{k} \mathcal{P}(g_{k+1}),Z_{n_{k+1}}=w_{k+1}
\end{array}
\right] \cdot  \bigl(1-\xi_{\tau(w_{k+1})}\bigr). 
\end{eqnarray*}
Thus we obtain the conditional probabilities:
\begin{eqnarray*}
& & \mathbb{P}_o\bigl[W_{k+1}=w_{k+1},\mathbf{i}_{k+1}=n_{k+1}\ |\
W_1^{k}=w_1^{k},\mathbf{i}_1^{k}=n_1^{k}\bigr] \\[1ex]
& = & \frac{1-\xi_{\tau(w_{k+1})}}{1-\xi_{\tau(w_{k})}} \cdot
\mathbb{P}_{w_k}\left[
\begin{array}{c}
\forall j\in \{1,\dots,n_{k+1}-2\}:Z_j\in C_{w_k},\\
Z_{n_{k+1}-1}\in
w_{k} \mathcal{P}(g_{k+1}),Z_{n_{k+1}}=w_{k+1}
\end{array}
\right] \\
& = & \frac{1-\xi_{\tau(w_{k+1})}}{1-\xi_{\tau(w_{k})}} \cdot \sum_{s\in \mathcal{P}(y)} \Bigl[
k_{\tau(w_k)}^{(n_{k+1}-1)}(o,s)\cdot p(s,y) \Bigr].
\end{eqnarray*}
\end{proof} 
Observe that the transition probabilities of the
stochastic process $(W_k,\mathbf{i}_k)_{k\in\mathbb{N}}$ depend only on 
 $\tau_k:=\tau(W_k)$, $\tau(W_{k+1})$ and $\widetilde W_{k+1}$. Hence, the stochastic process $\bigl(\widetilde
W_k,\mathbf{i}_k,\tau_k\bigr)_{k\in\mathbb{N}}$ is an irreducible Markov chain on the state space 
\begin{displaymath}
\mathcal{A}:=\Bigl\lbrace (y,n,j)\ \bigl|\ j\in\mathcal{I}, y\in V_j^\times,
n\in\mathbb{N}, \exists i\in\mathcal{I}\setminus\{j\}\exists s\in\mathcal{P}(y):k_i^{(n-1)}(o,s)>0\Bigr\rbrace
\end{displaymath}
with transition probabilities
\begin{displaymath}
q\bigl((x,m,i),(y,n,j)\bigr) =
\begin{cases}
0, & \textrm{if } i=j \\
\frac{1-\xi_j}{1-\xi_i}\cdot \sum_{s\in \mathcal{P}(y)}
\Bigl[ k_i^{(n-1)}(o,s)\cdot p(s,y)\Bigr], & \textrm{if } i\neq j
\end{cases}\,.
\end{displaymath}
For convenience, we write $q\bigl((x,m,i),(y,n,j)\bigr):=0$, if
$j\in\mathcal{I}$, $y\in V_j$, but $(y,n,j)\notin\mathcal{A}$. As the
probabilities $q\bigl((x,m,i),(y,n,j)\bigr)$ do not depend on $x$ and $m$, the sequence $(\tau_k)_{k\in\mathbb{N}}$ is
also a Markov chain on the state space $\mathcal{I}$ with transition probabilities
\begin{equation}\label{hatqij}
\hat q(i,j) := \sum_{y\in V_j^\times} \sum_{n\geq 1} q\bigl((x,m,i),(y,n,j)\bigr)
\end{equation}
for $i,j\in\mathcal{I}$ with $i\neq j$ and $\hat q(i,i)=0$. Note that $x\in V_i$ and $m\in\mathbb{N}$ can be
chosen arbitrarily, such that $(x,m,i)\in\mathcal{A}$. As $\mathcal{I}$ is finite, $(\tau_k)_{k\in\mathbb{N}}$ possesses an invariant probability
measure $\nu: \mathcal{I}\to [0,1]$, that is, for every $j\in\mathcal{I}$
\begin{equation}\label{nu}
\sum_{i\in\mathcal{I}} \nu(i) \cdot \hat q(i,j) = \nu(j)\,.
\end{equation}
We now define for $j\in\mathcal{I}$, $y\in V_j^\times$ and $n\in\mathbb{N}$
\begin{displaymath}
\pi(y,n,j) := \sum_{i\in\mathcal{I}} \nu(i)\cdot q\bigl((x,m,i),(y,n,j)\bigr),
\end{displaymath}
which is an invariant probability
measure of the stochastic process $\bigl(\widetilde W_k,\mathbf{i}_k,\tau_k\bigr)_{k\in\mathbb{N}}$, that is,
$$
\sum_{(x,m,i)\in\mathcal{A}} \pi(x,m,i)\, q\bigl((x,m,i),(y,n,j)\bigr) =
\pi(y,n,j)
$$
holds for all $(y,n,j)\in\mathcal{A}$.
\begin{Prop}\label{lambda}
There is a number $\Lambda\in \mathbb{R}$, such that 
\begin{displaymath}
\frac{\mathbf{e}_k}{k}\quad \xrightarrow{k\to\infty}\quad \Lambda \quad \mathbb{P}_o-a.s. .
\end{displaymath}
\end{Prop}
\begin{proof}
Consider the function $g: \mathcal{A} \rightarrow \mathbb{N}$, $(y,n,j)\mapsto
n$. An application of the ergodic theorem for positive recurrent Markov chains
shows that
\begin{displaymath}
\frac{1}{k}\sum_{l=1}^k g\bigl(\widetilde W_l,\mathbf{i}_l,\tau_l\bigr) =
\frac{\mathbf{e}_k-\mathbf{e}_0}{k}=\frac{\mathbf{e}_k}{k} \quad \xrightarrow{k\to\infty} \quad \int
g\,d\pi \quad \mathbb{P}_o-a.s.\,,
\end{displaymath}
if $\int g\,d\pi <\infty$ holds. Hence it is sufficient to show finiteness of
this integral. Noting that
\begin{eqnarray*}
\int g\,d\pi & = & \sum_{i\in\mathcal{I}} \nu(i) \cdot \sum_{(y,n,j)\in \mathcal{A}} n \cdot
q\bigl((x,m,i),(y,n,j)\bigr) \\
& = &  \sum_{i\in\mathcal{I}} \frac{\nu(i)}{1-\xi_i} \sum_{j\in\mathcal{I}\setminus\{i\}}
\alpha_j (1-\xi_j) \underbrace{\sum_{n\geq 1} n \sum_{y\in V_j^\times} \sum_{s\in \mathcal{P}(y)} k_i^{(n-1)}(o,s)\cdot p_j(s,y)}_{(*)},
\end{eqnarray*}
we now interpret the sum $(*)$ as a power series evaluated
at $1$. We have
\begin{eqnarray*}
&&\sum_{n\geq 1} n \sum_{y\in V_j^\times} \sum_{s\in \mathcal{P}(y)} k_i^{(n-1)}(o,s)\,
p_j(s,y)\,z^{n-1} \\
&=& 
\frac{\partial}{\partial z} \biggl[ \underbrace{\sum_{n\geq 1}  \sum_{y\in V_j^\times} \sum_{s\in \mathcal{P}(y)} k_i^{(n-1)}(o,s)\,
p_j(s,y)\,z^n}_{=: \gamma_{i,j}(z)} \biggr] .
\end{eqnarray*}
Now it is sufficient to show that the sum $\gamma_{i,j}(z)$ has radius of
convergence $R_{i,j}>1$ for all $i,j\in\mathcal{I}$ with $i\neq j$:
\begin{eqnarray*}
\gamma_{i,j}(z) &=& \sum_{n\geq 1} \sum_{y\in V_j^\times} \sum_{s\in \mathcal{P}(y)} k_i^{(n-1)}(o,s)\cdot p_j(s,y)\cdot
z^n\\
& = & \sum_{n\geq 1}  \sum_{s\in V_j} k_i^{(n-1)}(o,s)\cdot z^n - \sum_{n\geq 1}
\sum_{s\in \mathcal{P}(o_j)} k_i^{(n-1)}(o,s)\cdot p_j(s,o_j)\cdot z^n \\
& = & z\cdot \underbrace{\sum_{y\in V_j} K_i(o,y|z)}_{(**)} - z \underbrace{\sum_{s\in \mathcal{P}(o_j)} K_i(o,s|z)\cdot p_j(s,o_j)}_{(***)}. 
\end{eqnarray*}
From Equation (\ref{ki}) we obtain
$$
\sum_{y\in V_j} K_i(o,y|z) = \sum_{y\in V_j} \frac{L(o,y|z)}{1-\bar H_i(z)} =
\frac{1}{1-\bar H_i(z)} \sum_{y\in V_j} L_j\bigl(o_j,y|\xi_j(z)\bigr),
$$
and also
\begin{eqnarray*}
\sum_{y\in V_j} L_j\bigl(o_j,y|\xi_j(z)\bigr) & = & \frac{1}{G_j\bigl(o_j,o_j|\xi_j(z)\bigr)}
\sum_{y\in V_j} G_j\bigl(o_j,y|\xi_j(z)\bigr) \\
& = & \frac{1}{G_j\bigl(o_j,o_j|\xi_j(z)\bigr)} \cdot \frac{1}{1-\xi_j(z)}.
\end{eqnarray*}
Thus the sum $(**)$ has radius of convergence greater than $1$. Furthermore, by
(\ref{ki}) and Lemmas \ref{lemma-props} and \ref{lemma-fi},
\begin{eqnarray*}
&& \sum_{s\in \mathcal{P}(o_j)}p_j(s,o_j)\, K_i(o,s|z) \\
&=& \frac{1}{1-\bar H_i(z)} \sum_{s\in \mathcal{S}(o_j)}p_j(s,o_j)\,
\frac{G_j\bigl(o_j,s|\xi_j(z)\bigr)}{G_j\bigl(o_j,o_j|\xi_j(z)\bigr)} \\
&=& \frac{1}{\bigl(1-\bar H_i(z)\bigr)\cdot G_j\bigl(o_j,o_j|\xi_j(z)\bigr)} \sum_{s\in
    \mathcal{P}(o_j)}p_j(s,o_j)\, \xi_j(z)\, G_j\bigl(o_j,s|\xi_j(z)\bigr)\,
  \frac{1}{\xi_j(z)} \\
&=& \frac{G_j\bigl(o_j,o_j|\xi_j(z)\bigr)-1}{\bigl(1-\bar H_i(z)\bigr)\cdot G_j\bigl(o_j,o_j|\xi_j(z)\bigr)\cdot \xi_j(z)}\,.
\end{eqnarray*}
Thus the sum $(***)$ has radius of convergence greater than $1$, and also
$\gamma_{i,j}(z)$, whence $\int g\,d\pi$ is finite. 
\end{proof} 
Using the above, we can rewrite $\gamma_{i,j}(z)$ as
\begin{equation}\label{gammaij}
\gamma_{i,j}(z)  
 = \frac{1}{\alpha_i}\cdot 
  \frac{\xi_i(z)}{\xi_j(z)} \cdot 
  \biggl(\frac{1}{(1-\xi_j(z))\cdot G_j\bigl(o_j,o_j|\xi_j(z)\bigr)}
- 1 \biggr)
\end{equation}
and therefore
\begin{equation}\label{BigLambda}
\Lambda 
 =   \sum_{\substack{i,j\in\mathcal{I},\\ i\neq j}}  \nu(i)\cdot \alpha_j\cdot
\frac{1-\xi_j}{1-\xi_i}\cdot \gamma_{i,j}'(1).
\end{equation}
The following theorem is now obtained precisely as in Nagnibeda and Woess
\cite[proof of Theorem D]{woess2}:
\begin{Th}
\label{ell_fet}
\begin{displaymath}
\frac{\ell(Z_n)}{n} \quad 
\xrightarrow{n\to\infty} \quad \ell = \frac{1}{\Lambda} \quad \mathbb{P}_o-a.s.
\end{displaymath}
\end{Th}

Finally, we show how to compute the invariant probability measure
$\nu$ explicitly. For this purpose, it is sufficient to compute the transition probabilities $\hat
q(i,j)$ for all $i,j\in\mathcal{I}$. By solving the system of linear
equations given by Equations (\ref{nu}) $\nu$ is obtained. The next lemma proposes a formula for $\hat q(i,j)$:
\begin{Lemma}
Let $i,j\in\mathcal{I}$ with $i\neq j$. Then $\hat q(i,i)=0$ and
$$
\hat q(i,j) = \frac{\alpha_j}{\alpha_i} \cdot \frac{\xi_i}{\xi_j} \cdot
\frac{1-\xi_j}{1-\xi_i} \cdot \biggl( \frac{1}{(1-\xi_j) G_j(\xi_j)}-1
\biggr),\quad \textrm{ where } \ G_j(\xi_j):=G_j(o_j,o_j|\xi_j).
$$
\end{Lemma}
\begin{proof}
By definition of $W_k$ and $\mathbb{P}_o[Z_\infty\in V_\infty]=1$ it follows that $\hat
q(i,i)=0$. Considering (\ref{hatqij}) and the computations in the proof of Proposition \ref{lambda}
leads to $\hat q(i,j)= (1-\xi_j)/(1-\xi_i) \cdot \alpha_j\cdot \gamma_{i,j}(1)$.
\end{proof} 
We now give an explicit formula for $\nu$:
\begin{equation}\label{nu-formula}
\nu(i)=c\cdot \frac{\alpha_i\,(1-\xi_i)}{\xi_i} \cdot \bigl(1-(1-\xi_i)\,G_i(\xi_i)\bigr),
\end{equation}
where $c>0$ is chosen so that $\sum_{i\in\mathcal{I}} \nu(i)=1$. This is indeed
an invariant measure, because, writing $x(i)=1-(1-\xi_i)\,G_i(\xi_i)$, the
invariance condition on $\nu$ is just
$$
\sum_{i\in\mathcal{I}\setminus\{j\}} x(i) = \frac{x(j)}{\frac{1}{(1-\xi_j)
  G_j(\xi_j)} -1} \quad \textrm{ for each } j\in\mathcal{I}
$$
or, equivalently, that
$$
\sum_{i\in\mathcal{I}} x(i)=1.
$$
The following lemma verifies that this equation holds.
\begin{Lemma}\label{lemma-rho}
Let $i\in\mathcal{I}$. Then
$$
\rho(i):=\mathbb{P}_o\Bigl[ Z_\infty^{(1)}\in V_i^\times, \forall
n\in\mathbb{N}: Z_n\notin \bigcup_{j\in\mathcal{I}\setminus\{i\}} V_j\Bigr] =
\frac{1-(1-\xi_i)\, G_i(\xi_i)}{G(o,o|1)}.
$$
\end{Lemma}
\begin{proof}
By transience, $o$ is visited only finitely often $\mathbb{P}_e$-a.s., that is,
$$
\sum_{i\in\mathcal{I}} G(o,o|1)\,\rho(i)=1.
$$ 
This yields
\begin{eqnarray*}
\rho(i) & = & \sum_{y\in V_i^\times} \sum_{j\in\mathcal{I}\setminus\{i\}} L(o,y|1) \cdot \rho(j) \\
& = & \sum_{y\in V_i^\times}  L(o,y|1)\cdot  \bigl(G(o,o|1)^{-1}-\rho(i)\bigr) \\
& = & \bigl(G(o,o|1)^{-1}-\rho(i)\bigr) \sum_{y\in V_i^\times}
\frac{G_i\bigl(o_i,y|\xi_i\bigr)}{G_i\bigl(\xi_i\bigr)} \\
& = & \bigl(G(o,o|1)^{-1}-\rho(i)\bigr) \cdot \biggl( \frac{1}{(1-\xi_i) \cdot
  G_i(\xi_i)}-1\biggr).
\end{eqnarray*}
This leads to the proposed equation.
\end{proof}

We summarize this section: if we know for each factor $V_i$, $i\in\mathcal{I}$, the first visit generating function $F_i(x,o_i|z)$,
when starting at a predecessor $x\in\mathcal{P}(o_i)$ of $o_i$, and the Green function \mbox{$G_i(o_i,o_i|z)$},
it is possible to compute $\xi_i(z)$ by solving a finite system of characteristic
equations. This is in fact only possible when the generating functions are
known and not too complicated. The measure $\nu$ can then be computed by
(\ref{nu-formula}) and then $\ell$ found using $\ell=1/\Lambda$ and (\ref{BigLambda}). Sample computations are presented in
Section \ref{samples}.

%
%

\section{Double Generating Functions}
\label{bl_dgf}

In this section we compute the rate of escape $\ell$ for the random walk on $V$
 using double generating functions. The main tool
for our computation is the following theorem:

\begin{Th}[Sawyer and Steger]
\label{sawyer}
Let $(Y_n)$ be a sequence of real-valued random variables such that for some
$\delta>0$,
$$
\mathbb{E}\biggl( \sum_{n\geq 0} \exp(-rY_n-sn)\biggl) = \frac{C(r,s)}{g(r,s)}
\quad \textrm{ for } 0<r,s<\delta,
$$
where $C(r,s)$ and $g(r,s)$ are analytic for $|r|,|s|<\delta$
and \mbox{$C(0,0)\neq 0$}. Denote by $g_r$ and $g_s$ the partial derivatives of $g$
with respect to $r$ and $s$. Then
$$
\frac{Y_n}{n} \xrightarrow{n\to\infty} \ell=\frac{g_r(0,0)}{g_s(0,0)} \quad
\textrm{ almost surely.}
$$
\end{Th}
For the proof, see Sawyer and Steger \cite[Theorem 2.2]{sawyer}. 
\par
Setting $Y_n=\ell(Z_n)$, $w=e^{-r}$ and $z=e^{-s}$, to find $\ell$ in our
context it is sufficient to investigate the double generating
function
$$
\mathcal{E}(w,z):= \sum_{x\in V} \sum_{n\geq 0} p^{(n)}(o,x) w^{\ell(x)} z^n =
\sum_{x\in V} G(o,x|z) w^{\ell(x)}
$$
and to apply Theorem \ref{sawyer}. To this end, introduce further double
generating functions. Write $V^\times := V\setminus\{o\}$ and define
$$
\mathcal{L}(w,z) := 1+\sum_{n\geq 1}\ \sum_{x=x_1\dots x_n\in V^\times}\ \prod_{j=1}^n w\,
L_{\tau(x_j)} \bigl(o_{\tau(x_j)},x_j\bigl|\xi_{\tau(x_j)}(z)\bigr)
$$ 
and for $i\in\mathcal{I}$
\begin{eqnarray*}
\mathcal{L}_i^+(w,z) & := & \sum_{x\in V_i^\times} L_{i} \bigl(o_i,x\bigl|\xi_i(z)\bigr)\, w\,,\\
\mathcal{L}_i(w,z) & := & \mathcal{L}_i^+(w,z) \biggl( 1 + \sum_{n\geq 2} 
\sum_{\substack{x_2,\dots, x_n\in V^\times,\\ \tau(x_2)\neq i}}\, \prod_{j=2}^n w\,
L_{\tau(x_j)} \bigl(o_{\tau(x_j)},x_j\bigl|\xi_{\tau(x_j)}(z)\bigr)\biggr)\,.
\end{eqnarray*}
Thus we have the equation
$$
\mathcal{L}(w,z) = 1 + \sum_{i\in\mathcal{I}} \mathcal{L}_i(w,z).
$$
If $0\leq w,z<1$, the convergence of the series $\mathcal{L}(w,z)$ follows by \mbox{$\mathcal{L}(w,z)\leq \mathcal{E}(w,z)$}. Hence, $\mathcal{L}_i^+(w,z)$ and also
$\mathcal{L}_i(w,z)$, $i\in\mathcal{I}$, converge if \mbox{$|w|,|z|<1$}. The next
lemma provides another representation of $\mathcal{L}(w,z)$.

\begin{Lemma}
\label{scriptl}
Let $w,z\in \mathbb{R}$ with $0\leq w,z<1$. Then:
$$
\mathcal{L}(w,z) = \frac{1}{1-\mathcal{L}^\ast(w,z)}, \quad \textrm{ where }\quad  \mathcal{L}^\ast(w,z)=\sum_{i\in\mathcal{I}}
  \frac{\mathcal{L}_i^+(w,z)}{1+\mathcal{L}_i^+(w,z)}\,.
$$
\end{Lemma}
\begin{proof} 
Let $w,z\in \mathbb{R}$ with $0<w,z<1$. First we have
$$
\mathcal{L}_i(w,z) = \mathcal{L}_i^+(w,z) \cdot \Bigl(1+\sum_{j\in\mathcal{I}\setminus\{i\}}
\mathcal{L}_j(w,z)\Bigr)\quad \textrm{ for all } i\in\mathcal{I},
$$
and from convergence of $\mathcal{L}(w,z)$ we get
$$
\mathcal{L}_i(w,z) = \mathcal{L}_i^+(w,z) \cdot \bigl(\mathcal{L}(w,z) -
\mathcal{L}_i(w,z)\bigr)\,.
$$
As $\mathcal{L}_i^+(w,z)\geq 0$ holds, the last equation is equivalent to
$$
\mathcal{L}_i(w,z) =
\frac{\mathcal{L}_i^+(w,z)}{1+\mathcal{L}_i^+(w,z)} \mathcal{L}(w,z)\,.
$$
Thus
$$
\mathcal{L}(w,z)  =  1 + \sum_{i\in\mathcal{I}} \mathcal{L}_i(w,z)
 =  1 + \sum_{i\in\mathcal{I}} \frac{\mathcal{L}_i^+(w,z)}{1+ \mathcal{L}_i^+(w,z)} \mathcal{L}(w,z).
$$
As $\mathcal{L}(w,z)<\infty$, we get
$$
 \sum_{i\in\mathcal{I}} \frac{ \mathcal{L}_i^+(w,z)}{1+ \mathcal{L}_i^+(w,z)} < 1\,.
$$
and the result follows. 
\end{proof} 

\begin{Cor}
\label{scripte}
Let $w,z\in \mathbb{R}$ with $0\leq w,z<1$. Then
$$
\mathcal{E}(w,z)= \frac{G(o,o|z)}{1-\mathcal{L}^\ast(w,z)}\,.
$$
\end{Cor}
\begin{proof}
Let $w,z\in \mathbb{R}$ with $0<w,z<1$. Applying Lemma \ref{scriptl} yields the proposed equation:
\begin{eqnarray*}
\mathcal{E}(w,z) & = &  \sum_{x\in V} G(o,o|z) L(o,x|z) w^{\ell(x)} \\
& = & G(o,o|z) \Bigl( 1+ \sum_{n\geq 1} \sum_{x=x_1\dots x_n\in V^\times}
L(o,x|z) w^{\ell(x)} \Bigr) \\
& = & G(o,o|z) \cdot \mathcal{L}(w,z)\,.
\end{eqnarray*}
\end{proof} 
We can now conclude and compute a formula for the rate of escape
$\ell$. Rewriting $\mathcal{L}^\ast(w,z)$ after some manipulations
involving Lemma \ref{lemma-props} yields 
$$
\mathcal{L}^\ast(w,z)  = 
 \sum_{i\in\mathcal{I}} \frac{w \Bigl(\frac{1}{1-\xi_i(z)} - 
G_{i} \bigl(o_i,o_i|\xi_i(z)\bigr)\Bigr)}{\frac{w}{1-\xi_i(z)} +(1-w) G_{i} \bigl(o_i,o_i|\xi_i(z)\bigr)}\,.
$$
Now we define
\begin{eqnarray*}
C(w,z) & := & G(o,o|z) \quad \textrm{ and } \\
g(w,z) & := &  1-\sum_{i\in\mathcal{I}} \frac{w \Bigl(\frac{1}{1-\xi_i(z)} - 
G_{i} \bigl(o_i,o_i|\xi_i(z)\bigr)\Bigr)}{\frac{w}{1-\xi_i(z)} +(1-w) G_{i}
\bigl(o_i,o_i|\xi_i(z)\bigr)}\,.
\end{eqnarray*}
and we have
$$
\mathcal{E}(w,z) = \frac{C(w,z)}{g(w,z)} \quad \textrm{ for } 0\leq w,z<1.
$$
The constraints required for an application of Theorem \ref{sawyer} are obviously
fulfilled, as $G(o,o|z)$ has radius of convergence greater than $1$ and
$\xi_i<1$. We apply it now, where $g_w$ and $g_z$ denote the partial derivatives of
$g$ with respect to $w$ and $z$, respectively. Hence, we can conclude
$$
\frac{\ell(Z_n)}{n} \xrightarrow{n\to\infty} \ell=\frac{g_w(1,1)}{g_z(1,1)}
\quad \mathbb{P}_o-a.s.\,.
$$
Simplifications yield the following formula for $\ell$, where we write
$\xi_i=\xi_i(1)$ and $G_i(\xi_i):=G_i(o_i,o_i|\xi_i)$:
\begin{equation}
\label{roe-dgf}
\ell = \frac{\sum_{i\in\mathcal{I}} \Bigl[ \bigl(1-(1-\xi_i)\,
  G_i(\xi_i)\bigr)\cdot  G_i(\xi_i)\cdot (1-\xi_i) \Bigr]}{
  \sum_{i\in\mathcal{I}} \Bigr[ \xi_i'(1)\cdot \bigl( G_i(\xi_i) -(1-\xi_i)\, G_i'(\xi_i)\bigr) \Bigr] }\,.
\end{equation}
Observe that $\xi_i,G_i(\xi_i)>0$, $\xi_i<1$ and $G_i(\xi_i)<(1-\xi_i)^{-1}$. Thus $\ell>0$.

%
%

\section{Free Products of Groups}

\label{bl_lp}

In this section we present a third technique for the computation of the rate of
escape $\ell$ of the block length for the random walk on the free product. This
technique is restricted to the case of a free product of \textit{groups}. Therefore let
$\Gamma_i$, $i\in\mathcal{I}$, be non-trivial finitely generated groups. We assume that
the groups have pairwise trivial intersections, but they may be isomorphic. Denote by $e_i$ the identity on
$\Gamma_i$. The elements of the free product $\Gamma:=\Gamma_1\ast \dots\ast
\Gamma_r$ are represented as words in the sense of (\ref{freeproduct}) and $e$
is identified with the empty word. 
\par
We can define a group operation on $\Gamma$: the product of $u,v\in\Gamma$ is the concatenation of the words $u$
and $v$ with possible cancellations and contractions in the middle to get the
representative form of the product word. We exclude the case
$r=2=|\Gamma_1|=|\Gamma_2|$. This ensures that the free group product is
non-amenable, yielding that each of our constructed random walks on $\Gamma$ is transient and $G(e,e|z)$ has
radius of convergence greater than $1$. (See Woess \cite[Theorem 10.10,
Proposition 12.4, Corollary 12.5]{woess}).
We write $\Gamma_i^\times
=\Gamma_i\setminus \{e_i\}$ for $i\in\mathcal{I}$. Furthermore we write
$\Gamma_\infty = V_\infty$ and $\Gamma^\times:=\Gamma \setminus\{e\}$. 
\par
The random walk on $\Gamma$ is constructed as follows: standing at $x\in\Gamma$
we allow walking to $xg$ with $g\in\bigcup_{i\in\mathcal{I}}\Gamma_i^\times$ in one step. Choose probability measures
$\mu_i$ on $\Gamma_i^\times$ for each $i\in\mathcal{I}$ such that $\mu_i$ defines an
irreducible random walk on $\Gamma_i$, that is, $p_i(x,y)=\mu_i(x^{-1}y)$ for
all $x,y\in\Gamma_i$. Let  $\alpha_1,\dots,\alpha_r>0$ with $\sum_{i\in\mathcal{I}}
\alpha_i=1$. Then we define the transition probabilities as
$$
p(x,xg) := \alpha_i \cdot \mu_i(g)
$$
for all $x\in \Gamma$, $g\in \Gamma_i^\times$, and we set $p(x,y):=0$ otherwise. As the transition
probabilities depend only on the increment $g\in \Gamma_i^\times$ we write $\mu(g):=p(x,xg)$ for all
\mbox{$g\in \bigcup_{i\in\mathcal{I}} \Gamma_i^\times$}, and $\mu(g):=0$
otherwise. Analogously, the n-step transition probabilities are given by the
convolution powers $\mu^{(n)}$ of $\mu$. 
\par
By Theorem
\ref{elltoinfty} the random walk converges $\mathbb{P}_e$-a.s. to a random
variable $Z_\infty$ with values in $\Gamma_\infty$. Denote by $\nu$ the
distribution of $Z_\infty$. Let 
$$
E_i=\bigl\lbrace x_1x_2\dots \in \Gamma_\infty\ \bigl|\ \tau(x_1)=i\bigr\rbrace \, \textrm{ for } i\in\mathcal{I}.
$$
Then $\nu$ is uniquely determined by its values on the Borel sets $B$ of the form $xE_i=\{xh \mid h\in E_i\}$ with $i\in\mathcal{I}$,
$x\in\Gamma$ and $\tau(x)\neq i$. We will now give a formula for these values:

\begin{Lemma}
Let $i\in\mathcal{I}$, $x\in \Gamma$ with $\tau(x)\neq i$. Then
$$
\nu(xE_i) = \mathbb{P}_e[Z_\infty \in xE_i] = F(e,x|1) \bigl(1-(1-\xi_i)\, G_i(e_i,e_i|\xi_i)\bigr)\,.
$$
\end{Lemma}
\begin{proof}
The proof of this lemma is extrapolated from Woess \cite[Theorem 4\,c]{woess1},
where one can find an incorrect formula, which we correct here. First we have
$$\nu(xE_i) = F(e,x|1) \cdot \nu(E_i).
$$
Recall that we have by vertex-transitivity $G_i(o_i,o_i|z)=G_i(y,y|z)$
for all $i\in\mathcal{I}$ and all $y\in\Gamma_i$. By Lemma \ref{lemma-rho} we obtain
$$
\nu(E_i) = G(o,o|1) \cdot \rho(i) = 1-(1-\xi_i)\, G_i(e_i,e_i|\xi_i).
$$
This leads to the proposed formula. 
\end{proof} 
Now we reformulate our problem for finding a formula for $\ell$. For this
purpose, we apply a technique which was used by Ledrappier
\cite[Section 4\,b]{ledrappier} for free groups. 
\par
By Corollary
\ref{ell_fet} and Lebesgue's Dominated Convergence Theorem we have $\mathbb{P}_e-a.s.$
$$
\lim_{n\to\infty} \frac{\mathbb{E}[\ell(Z_n)]}{n} = \lim_{n\to\infty} \int
\frac{\ell(Z_n)}{n} \,d\mathbb{P}_e = \int \lim_{n\to\infty}
\frac{\ell(Z_n)}{n} d\mathbb{P}_e = \int \ell\, d\mathbb{P}_e = \ell .
$$
Thus it is sufficient to prove convergence of the sequence 
$$
\Bigl( \mathbb{E}[\ell(Z_{n+1})] - \mathbb{E}[\ell(Z_{n})] \Bigr)_{n\in\mathbb{N}}
$$
and to compute its limit, which then must equal $\ell$. First we have
$$
\mathbb{E}[\ell(Z_{n})] = \sum_{h\in \Gamma} \ell(h)\, \mu^{(n)}(h)
$$
and 
$$
\mathbb{E}[\ell(Z_{n+1})] = \sum_{g,h\in \Gamma} \ell(gh)\,
\mu(g)\, \mu^{(n)}(h).
$$
On the other hand,
$$
\mathbb{E}[\ell(Z_{n})] = \sum_{g\in \Gamma} \mu(g)\,
\mathbb{E}[\ell(Z_{n})] = \sum_{g,h\in \Gamma} \mu(g)\,
\ell(h)\, \mu^{(n)}(h) .
$$
Thus we obtain
\begin{eqnarray*}
\mathbb{E}[\ell(Z_{n+1})] -\mathbb{E}[\ell(Z_{n})] & = & \sum_{g\in \Gamma} \mu(g) \sum_{h\in \Gamma} \bigl(\ell(gh)-\ell(h)\bigr)\, \mu^{(n)}(h) \\
& = &  \sum_{g\in \Gamma} \mu(g) \int_\Gamma \bigl( \ell(gZ_n)
-\ell(Z_n)\bigr)\, d\mathbb{P}_e .
\end{eqnarray*}
Define now the random variables
$$
Y_n := \ell(gZ_n) -\ell(Z_n)
$$
for any given $g\in\bigcup_{i\in\mathcal{I}} \Gamma_i^\times$. We have $Y_n\in\{-1,0,1\}$ for all
$n\in\mathbb{N}$.
By vertex-transitivity $gZ_n$ converges to $gZ_\infty$. Hence, $Y_n$ converges to a
random variable $Y_\infty$ with values in $\{-1,0,1\}$ depending only on $g$ and
the first block of $Z_\infty$. In other words, $Y_n$ becomes constant, if
$n$ is big enough. If $Z_\infty=x_1x_2\dots$, we obtain
for given $g\in\bigcup_{i\in\mathcal{I}} \Gamma_i^\times$:
$$
Y_\infty = 
\begin{cases}
0 & \textrm{, if } \tau(x_1)=\tau(g)\ \textrm{ and }\ x_1g\neq e \\
-1 & \textrm{, if }\tau(x_1)=\tau(g)\ \textrm{ and }\ x_1g= e \\
1 & \textrm{, if }\tau(x_1)\neq \tau(g)
\end{cases}.
$$
By Lebesgue's Dominated Convergence Theorem, we infer that
$$
\int \bigl( \ell(gZ_n) -\ell(Z_n)\bigr) \,d\mathbb{P}_e \xrightarrow{n\to\infty} \int
Y_\infty \, d\mathbb{P}_e.
$$
Consider the function
$$
f: \biggl(\bigcup_{i\in\mathcal{I}} \Gamma_i^\times \biggr) \times \Gamma_\infty \rightarrow \{-1,0,1\}:
(g,x_1x_2\dots )
\mapsto  \ell(gx_1)-\ell(x_1) 
$$
and its projections $f_g: \Gamma_\infty \rightarrow \{-1,0,1\}: w \mapsto
f(g,w)$ for every $g\in\bigcup_{i\in\mathcal{I}} \Gamma_i^\times$. Observe that each $f_g$ is measurable and thus
$$
\int_{\Gamma_\infty} Y_\infty\, d\mathbb{P}_e = \int_{\Gamma_\infty} f(g,Z_\infty)\,
d\nu = \int_{\Gamma_\infty} f_g(w)\, d\nu(w)\,.
$$
Denote by $E_{h}$ the event that $Z_\infty$ has as first block the element
$h\in\bigcup_{i\in\mathcal{I}} \Gamma_i^\times$ and denote $E_{\neq i}$ the event that
$Z_\infty$ starts with a block element not of type $i\in\mathcal{I}$. Then we
obtain for $g\in\Gamma_i$
$$
\nu(E_{g^{-1}}) = F(e,g^{-1}|1) \cdot \bigl(1-\nu(E_i)\bigr) =  F(e,g^{-1}|1)
\cdot (1-\xi_i)\cdot G_i(e_i,e_i|\xi_i)
$$
and
$$
\nu(E_{\neq i})  =  (1-\xi_i)\cdot G_i(e_i,e_i|\xi_i)\bigr)\,.
$$
Writing $G_i(z)$ for
$G_i(e_i,e_i|z)$, as before,  we can conclude:
\begin{eqnarray*}
& & \mathbb{E}[\ell(Z_{n+1})] -\mathbb{E}[\ell(Z_{n})] \\
 & \xrightarrow{n\to\infty} & \sum_{g\in \bigcup_{i\in\mathcal{I}} \Gamma_i^\times} \mu(g)
 \int_{\Gamma^\infty} f_g(w)\, d\nu(w) \\
& = & \sum_{g\in \bigcup_{i\in\mathcal{I}} \Gamma_i^\times} \mu(g) \bigl( - \nu(E_{g^{-1}}) +
\nu(E_{\neq i})\bigr) \\
& = & \sum_{i\in\mathcal{I}}\alpha_i\,(1-\xi_i)\,
G_i(\xi_i)\,\bigl(1-\underbrace{\sum_{g\in\Gamma_i^\times}
  \mu_i(g)\,F_i(e_i,g^{-1}|\xi_i)}_{=\frac{G_i(\xi_i)-1}{\xi_i\cdot
  G_i(\xi_i)}}\bigr)\\
&=& \sum_{i\in\mathcal{I}}\alpha_i\,\frac{1-\xi_i}{\xi_i}\,\bigl(1-(1-\xi_i)\,G_i(\xi_i)\bigr).
\end{eqnarray*}
Thus we get the rate of escape of the block length as
\begin{equation}\label{roe-lp}
\ell  =  \sum_{i\in\mathcal{I}}\alpha_i\,\frac{1-\xi_i}{\xi_i}\,\bigl(1-(1-\xi_i)\,G_i(\xi_i)\bigr).
\end{equation}
Observe that the technique presented in this section can be extended to a free
product of an infinite, countable number of groups. All the required properties
of the generating functions used also hold in this
case. Furthermore, $Y_{n}$ is again bounded such that finiteness of $\ell$ is
ensured. Thus the same computations prove the same formula for $\ell$ if $\mathcal{I}=\mathbb{N}$.

%
%

\section{Examples}
\label{samples}
We present two sample applications of our formulas for the rate of escape of
the block length. First we look at a free product arising from non-Cayley graphs and then
we look at a free product of infinite groups. Note
that both examples go beyond previously investigated graph structures for the
computation of $\ell$.

\subsection{Free Product arising from Non-Cayley-Graphs}
\label{noncayleysample}
\begin{figure}[htp]
\begin{center}
\includegraphics[width=7cm]{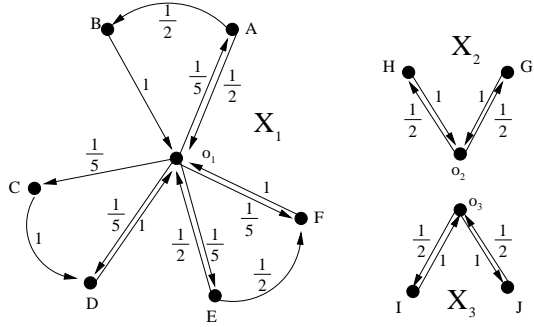}
\end{center}
\caption{Some non-Cayley graphs}
\label{graphsample}
\end{figure}

Consider the sets $V_1=\{A,B,C,D,E,F,o_1\},V_2=\{G,H,o_2\}$ and $V_3=\{I,J,o_3\}$ and
the random walks on these sets. Their transition probabilities are
sketched in Figure \ref{graphsample}. Note that none of the graphs in this figure is a Cayley graph. 
Consider now the corresponding random walk on the free product $V=V_1\ast
V_2\ast V_3$, where $\alpha_1=5/9$ and $\alpha_2=\alpha_3=2/9$. 
\par
We obtain the following generating functions:
\begin{eqnarray*}
U_1(o_1,o_1|z) & = & \frac{3}{5} z^2 +\frac{2}{5} z^3,\\
G_1(o_1,o_1|z) & = & \frac{1}{1-U_1(o_1,o_1|z)} = \frac{1}{1-\frac{3}{5}z^2
  -\frac{2}{5}z^3}, 
\end{eqnarray*} 
\begin{eqnarray*}
F_1(A,o_1|z) & = & F_1(E,o_1|z) = \frac{1}{2}z + \frac{1}{2} z^2, \\
F_1(C,o_1|z) & = & z^2, \quad F_1(D,o_1|z)= F_1(F,o_1|z)=z,\\
G_2(o_2,o_2|z) & = & \frac{1}{1-z^2},\\
F_2(G,o_2,z) & = & F_2(H,o_2|z) = z, \\
\bar H_1(z) & = &  \frac{4}{9}z \xi_2(z),\\
\bar H_2(z) & = & \bar H_3(z) = 
\frac{2}{9} z \xi_2(z)+ \frac{5}{9}\cdot \frac{1}{5} z \bigl(3\xi_1(z) + 2\xi_1(z)^2\bigr).
\end{eqnarray*} 
Note that $\bar H_2(z)=\bar H_3(z)$ follows by symmetry. The 
Green function $G(o,o|z)$ of the corresponding random walk on $V$ has radius of convergence greater than $1$. This can be
shown by constructing recursive equations using $\bar H_1(z)$ and $\bar H_2(z)$
and numerical evaluation. Consider
\begin{eqnarray*}
\xi_1(z) & = & \frac{\frac{5}{9} z}{1-\bar H_1(z)} = \frac{\frac{5}{9} z}{1-\frac{4}{9} z \xi_2(z)} \quad \textrm{ and} \\
\xi_2(z) & = & \xi_3(z) = \frac{\frac{2}{9} z}{1-\bar H_2(z)} = \frac{\frac{2}{9}
  z}{1-\frac{2}{9} z \xi_2(z) - \frac{1}{9} z \bigl(3\xi_1(z)+2\xi_1(z)^2\bigr)}.
\end{eqnarray*} 
Substituting $\xi_1(z)$ into $\xi_2(z)$ we have to solve an equation in
the variable $\xi_2(z)$. Solving this equation with \textsc{mathematica} we obtain four continuous solutions, but only
one solution satisfies $\xi_2(1)<1$. Hence,
we get $\xi_2(z)$ as this solution and obtain $\xi_1(z)$ from this
$\xi_2(z)$. We find that 
$$
\xi_1(1)\approx 0.66571\quad \textrm{ and } \quad \xi_2(1)=\xi_3(1)\approx 0.37231.
$$
We compute $\ell$ using Theorem \ref{ell_fet}. The transition matrix of the
Markov chain $(\tau_k)_{k\in\mathbb{N}}$ of the
alternating vertex types is
$$
\bigl(\hat q(i,j)\bigr)_{1\leq i,j \leq 3} = \left( 
\begin{array}{ccc}
0 & 0.5 & 0.5 \\
0.62769 & 0 & 0.37231 \\
0.62769 & 0.37231 & 0 
\end{array}
\right)
$$
and from this we obtain the corresponding invariant probability measure $\nu$ with
$$
\nu(1) = 0.38563 \quad \textrm{and} \quad \nu(2) = \nu(3) = 0.30718.
$$
Now we are able to compute the rate of escape of the block length to the
random walk on $V$. We obtain
$$
\ell \approx 0.33089.
$$
If $\ell$ is computed by Equation (\ref{roe-dgf}), then the numerical approximated result
and the above result coincide in the first 50 decimal numbers. So numerical
approximations do not lead to a distortion of the result.

\subsection{$\mathbb{Z}^2\ast \mathbb{Z}/2$}
Consider $V_1=\mathbb{Z}^2$ and the
simple random walk on it given by $\mu_1\bigl((\pm1,0)\bigr)=\mu_1\bigl((0,\pm 1)\bigr)=1/4$. Also,
consider the group $V_2=\mathbb{Z}/2$ and the simple random walk on it
given by $\mu_2(1_2)=1$. We are interested now in the simple random walk on
$V=V_1\ast V_2$, where $\alpha_1=4/5$ and $\alpha_2=1/5$. For the
computation of $\ell$ we use Equation (\ref{roe-lp}). Therefore it is sufficient to compute $\xi_1,\xi_2$ and
$G_1\bigl((0,0),(0,0)|\xi_1\bigr)$. For this purpose, we use the computations
and results in Woess \cite[pages 100, 105, 109]{woess}; compare also with
Soardi \cite{soardi}.
\par
\begin{figure}
\includegraphics[width=5cm]{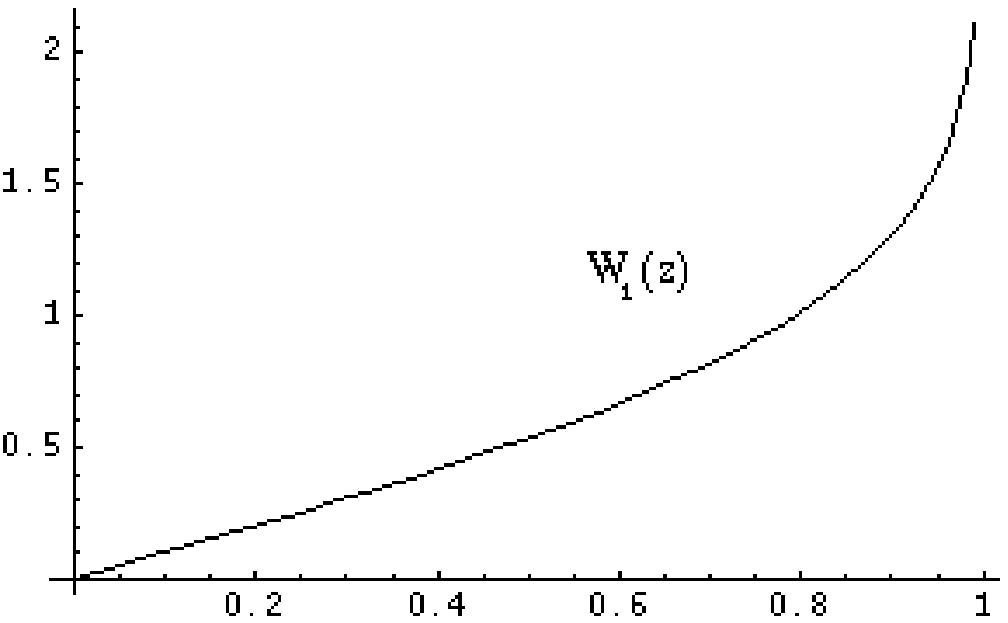}
\hspace{1cm}
\includegraphics[width=5cm]{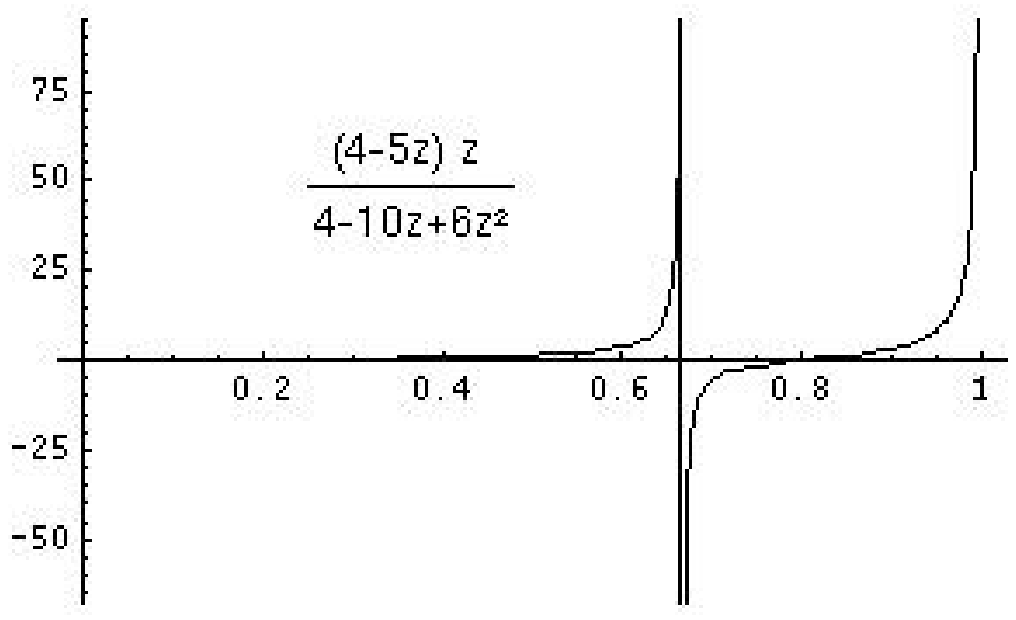}
\caption{Graphs of $W_1(z)$ and $(4-5z)z/(4-10z+6z^2)$}
\label{xi-graphs}
\end{figure}

Before we can compute these values, we have to introduce some auxiliary
functions. In the following let the subindex $0$ correspond to the random walk
on $V$. Denote 
$$
W_i(z)=z\cdot G_i(o,o|z) \textrm{
  for } i\in\{0,1,2\}.
$$
As
$W_i(z)$ is strictly increasing, there is an inverse function
$W^{-1}_i(z)$ such that \mbox{$W^{-1}_i(W_i(z))=z$} holds. By \cite[Theorem 9.10]{woess} we have 
$$
G_i(o,o|z)=\Phi_i\bigr(z G_i(o,o|z)\bigr),\quad \textrm{ where }
\Phi_i(t)=\frac{t}{W^{-1}_i(t)} \quad \forall  i\in\{0,1,2\}.
$$
By \cite[Example 9.15 (3)]{woess}, 
$$
W_1(z) = \frac{1}{4\pi^2} \int_{(-\pi,\pi]^2} \frac{2z}{2-z\cdot (\cos x_1
  +\cos x_2)} d\underline{x},
$$
where $\underline{x}=(x_1,x_2)$. Furthermore, by \cite[Theorem 9.19]{woess} we have the equation
$$
\Phi_0(t) = \Phi_1(\alpha_1 t) + \Phi_2(\alpha_2 t) -1.
$$
By \cite[Example 9.15 (1)]{woess} we have
$$
\Phi_2(t) = \frac{1}{2} \Bigl( \sqrt{1+4t^2} +1\Bigr).
$$
This yields
$$
\Phi_1\Bigl(\frac{4}{5}W(z)\Bigr) = G(o,o|z) - \frac{1}{2} \biggl( \sqrt{1+\frac{4}{25} W(z)^2} -1\biggr).
$$
Inverting this equation and multiplication with $4/5\cdot W(z)$ leads to
$$
\frac{\frac{4}{5} W(z)}{\Phi_{1} \bigl(\frac{4}{5}W(z)\bigr)} = \frac{\frac{4}{5}W(z)}{G(o,o|z) -\frac{1}{2} \Bigl( \sqrt{1+\frac{4}{25} W(z)^2} -1\Bigr)}.
$$
Applying $W_1$ onto both sides of this equation yields
$$
\frac{4}{5} W(z) = W_1\Biggl( \frac{\frac{4}{5}W(z)}{G(o,o|z) -\frac{1}{2}
  \Bigl( \sqrt{1+\frac{4}{25} W(z)^2} -1\Bigr)} \Biggr).
$$
Note that we have
\begin{equation}\label{w1-equ}
W_1\bigl(\xi_1(z)\bigr) = \frac{4}{5} W(z).
\end{equation}
and thus
$$
\xi_1(z) = \frac{\frac{4}{5} W(z)}{\frac{W(z)}{z} - \Bigl(\frac{1}{2}
  \sqrt{1+\frac{4}{25} W(z)^2}-1\Bigr)}.
$$
Substituting $y=4/5\cdot W(z)$ we obtain
$$
\xi_i=\xi_1(1) = \frac{y}{\frac{5y}{4} - \Bigl(\frac{1}{2}
  \sqrt{1+\frac{1}{4} y^2}-1\Bigr)}
$$
or equivalently
$$
y= \frac{(4 - 5\,\xi_1)\ \xi_1}{4 - 10\, \xi_1 + 6\, \xi_1^2}
$$
Substituting $y$ into Equation (\ref{w1-equ}), we have to solve
$$
\frac{(4 - 5\,\xi_1)\ \xi_1}{4 - 10\,\xi_1 + 6\,\xi_1^2} = W_1(\xi_1)
$$
in the unknown variable $\xi_1$. The solution can be computed
only numerically. Considering the graphs of the functions $(4-5z)z/(4-10z+6z^2)$
and $W_1(z)$ we see that there is only one possible intersection point greater
than $4/5$. See Figure \ref{xi-graphs}. 
Using the bisection method and numerical integration and evaluation we obtain
$$
\xi_1 \approx 0.84426, \textrm{ and also }
W(1)=\frac{5}{4} y \approx 1.40724.
$$
This yields
$$
G_1\bigl((0,0),(0,0)|\xi_1\bigr)= \frac{\frac{4}{5}W(1)}{\xi_1(1)} \approx 1.33347.
$$
By
$$
\xi_1= \frac{\frac{4}{5}}{1-\frac{1}{5}\xi_2}.
$$
we obtain
$$
\xi_2 \approx 0.26212, \textrm{ and }
G_2(0,0|\xi_2) = \frac{1}{1-\xi_2^2} \approx 1.07378.
$$
Now we have computed all necessary characteristical numbers and the rate of
escape of the block length of the simple random walk on $\mathbb{Z}^2\ast
\mathbb{Z}/2$ can be computed by Equation (\ref{roe-lp}) as
$$
\ell \approx 0.23386.
$$

\section{Remarks} \label{remarks}
We can extend our considerations in order to compute other types of typical
rates of escapes concerning the random walk on the
free product:
\begin{itemize}
\item[I.] For $x=x_1\dots x_n\in V$ and $i\in\mathcal{I}$ the \textsf{partial block
length of $x$ with respect to $V_i$} is given by
$$
\ell_i(x) := \bigl| \bigl\lbrace j\ |\ j\in\{1,\dots,n\}, x_j\in V_i\bigr\rbrace\bigr|.
$$
As $\ell_i(W_k)/k$ converges for $k\to\infty$ to $\nu(i)$, which is
the invariant probability measure on $\mathcal{I}$ with respect to the Markov chain
$(\tau(W_k))_{k\in\mathbb{N}}$, we obtain the \textsf{partial rate of escape of
  the block length}
$$
\frac{\ell_i(Z_n)}{n} \quad 
\xrightarrow{n\to\infty} \quad \nu(i) \cdot \ell \quad \mathbb{P}_o-a.s..
$$

\item[II.] 
The set $V$ carries a \textit{Markovian distance} defined by 
$$
d(x,y):=\min \bigl\lbrace n\in\mathbb{N}\ |\ p^{(n)}(x,y)>0\bigr\rbrace,
$$
where $x,y\in V$. Note that in general $d(\cdot,\cdot)$ is not necessarily
symmetric. The \textsf{Markovian length} is
defined as $|x|=d(o,x)$. We can extend the considerations of Section \ref{bl_fet} for the proof of the existence of the rate of escape of the Markovian length, that is
$$
\lambda = \lim_{n\to\infty} \frac{1}{n} |Z_n|.
$$
This yields also a formula for $\lambda$ given by
$$
\lambda = \ell \cdot \sigma, \textrm{ where } \sigma = \sum_{\substack{i,j\in\mathcal{I},\\ i\neq j}} \nu(i) 
 \frac{\alpha_j}{\alpha_i} \frac{\xi_i}{\xi_j} \frac{1-\xi_j}{1-\xi_i}
 \sum_{m\geq 1} m \cdot \sum_{y\in
 S_{j}(m)} L_j(o_j,y|\xi_j)
$$
and $S_{i}(m)=\bigl\lbrace x\in V_i\ \bigl|\ |x|=m\ \bigr\rbrace$ for
$i\in\mathcal{I}$ and $m\in\mathbb{N}$.
\end{itemize}

\section*{Acknowledgement}
The author is grateful to Wolfgang Woess for numerous discussions on several
problems and his help during the preparation of this article, and also to
Donald Cartwright for several hints regarding content and exposition.

\bibliographystyle{abbrv}
\bibliography{literatur}

\begin{thebibliography}{10}

\bibitem{cartwright-kaimanovich-woess}
D.~Cartwright, V.~Kaimanovich, and W.~Woess.
\newblock Random walks on the affine group of local fields and of homogenous
  trees.
\newblock {\em Ann. Inst. Fourier (Grenoble)}, 44:1243--1288, 1994.

\bibitem{cartwright-soardi}
D.~Cartwright and P.~Soardi.
\newblock Random walks on free products, quotients, and amalgams.
\newblock {\em Nagoya Math. J.}, 102:163--180, 1986.

\bibitem{derrienic}
Y.~Derriennic.
\newblock Quelques applications du th\'eor\`eme ergodique sous-additif.
\newblock {\em Ast\'erisque}, 74:183--201, 1980.

\bibitem{dyubina}
A.~Dyubina.
\newblock Characteristics of random walks on wreath products of groups.
\newblock {\em J. of. Math. Sciences}, 107(5):4166--4171, 2001.

\bibitem{erschler2}
A.~Erschler.
\newblock On the asymptotics of drift.
\newblock {\em J. of Math. Sciences}, 121(3):2437--2440, 2004.

\bibitem{furstenberg}
H.~Furstenberg.
\newblock Non commuting random products.
\newblock {\em Trans. Amer. Math. Soc.}, 108:377--428, 1963.

\bibitem{guivarch}
Y.~Guivarc'h.
\newblock Sur la loi des grands nombres et le rayon spectral d'une marche
  al\'eatoire.
\newblock {\em Ast\'erisque}, 74:47--98, 1980.

\bibitem{kaimanovich-vershik}
V.~Kaimanovich and A.~Vershik.
\newblock Random walks on discrete groups: boundary and entropy.
\newblock {\em Ann. Probab.}, 11:457--490, 1983.

\bibitem{kingman}
J.~Kingman.
\newblock The ergodic theory of subadditive processes.
\newblock {\em J. Royal Stat. Soc., Ser. B}, 30:499--510, 1968.

\bibitem{ledrappier}
F.~Ledrappier.
\newblock Some asymptotic properties of random walks on free groups.
\newblock In {\em CRM Proceedings and Lecture Notes}, volume~28, pages
  117--152. CRM, 2001.

\bibitem{lyons-pemantle-peres}
R.~Lyons, R.~Pemantle, and Y.~Peres.
\newblock Random walks on the lamplighter group.
\newblock {\em Ann. of Probability}, 24(4):1993--2006, 1996.

\bibitem{mairesse}
J.~Mairesse.
\newblock Random walks on groups and monoids with a markovian harmonic measure.
\newblock {\em Electron. J. Probab.}, 10:1417--1441, 2005.

\bibitem{mairesse-mantaray}
J.~Mairesse.
\newblock Randomly growing braid on three strands and the manta ray.
\newblock {\em Ann. Appl. Probab.}, 17(2):502--536, 2007.

\bibitem{mairesse1}
J.~Mairesse and F.~Math\'eus.
\newblock Random walks on free products of cyclic groups.
\newblock {\em J. London Math. Soc.}, 75(1):47--66, 2007.

\bibitem{mclaughlin}
J.~McLaughlin.
\newblock {\em Random walks and convolution operators on free products}.
\newblock PhD thesis, New York Univ., 1986.

\bibitem{woess2}
T.~Nagnibeda and W.~Woess.
\newblock Random walks on trees with finitely many cone types.
\newblock {\em J. Theoret. Probab.}, 15:399--438, 2002.

\bibitem{sawyer}
S.~Sawyer and T.~Steger.
\newblock The rate of escape for anisotropic random walks in a tree.
\newblock {\em Probab. Theory Related Fields}, 76:207--230, 1987.

\bibitem{soardi}
P.~Soardi.
\newblock Simple random walks on $\mathbb{Z}^2\ast \mathbb{Z}_2$.
\newblock {\em Symposia Math.}, 29:303--309, 1986.

\bibitem{varopoulos}
N.~T. Varopoulos.
\newblock Long range estimates for markov chains.
\newblock {\em Bull. Sc. math.}, 109:225--252, 1985.

\bibitem{voiculescu}
D.~Voiculescu.
\newblock Addition of certain non-commuting random variables.
\newblock {\em J. Funct. Anal.}, 66:323--346, 1986.

\bibitem{woess1}
W.~Woess.
\newblock A description of the martin boundary for nearest neighbour random
  walks on free products.
\newblock {\em Probability Measures on Groups}, VII:203--215, 1985.

\bibitem{woess3}
W.~Woess.
\newblock Nearest neighbour random walks on free products of discrete groups.
\newblock {\em Boll. Un. Mat. Ital.}, 5-B:961--982, 1986.

\bibitem{woess}
W.~Woess.
\newblock {\em Random Walks on Infinite Graphs and Groups}.
\newblock Cambridge University Press, 2000.

\end{thebibliography}

\end{document}